\documentclass[a4paper,11pt]{amsart}
\usepackage{mathrsfs}
\usepackage{bbm}
\usepackage{latexsym, amssymb,amsmath,amsthm}

\def \To{\longrightarrow}

\def \Hom{\operatorname{Hom}}
\def \Ext{\operatorname{Ext}}

\def \g{\mathfrak{g}}

\def \sgn{\operatorname{sgn}}
\def \H{\operatorname{H}}
\def \Ker{\operatorname{Ker}}
\def \Im{\operatorname{Im}}

\numberwithin{equation}{section}

\newtheorem{theorem}{Theorem}[section]
\newtheorem{lemma}[theorem]{Lemma}
\newtheorem{proposition}[theorem]{Proposition}

\newtheorem{definition}[theorem]{Definition}

\newtheorem{remark}[theorem]{Remark}

\begin{document}
\title[6-term exact sequence]
{A Hochschild's 6-term exact sequence for restricted Lie
superalgebras}
\author{Gongxiang Liu}
\address{Department of Mathematics, Nanjing University, Nanjing 210093, China} \email{gxliu@nju.edu.cn}
\maketitle

\begin{abstract} In \cite{Hoch}, Hochschild established a 6-term
exact sequence for the cohomology of restricted Lie algebras. We
generalize this result to restricted Lie superalgebras.
\end{abstract}

\section{Introduction}
\subsection{} As  generalizations and deep continuations of classical
Lie theory, Lie superalgebras over the field of complex numbers
$\mathbb{C}$ have been studied extensively since the classification
of finite dimensional complex simple Lie superalgebras by Kac
\cite{Kac1}. Comparing with abundant works and results for various
cohomology theory of Lie superalgebras (see \cite{BKN,HP,KK,SZ} and
references therein), the knowledge about the cohomology theory of
restricted Lie superalgebras is poor. To the author's best
knowledge, there has not been any serious study in the direction,
perhaps because even for simple Lie superalgebras over $\mathbb{C}$
their cohomology theory is already very difficult.

\subsection{} In \cite{Hoch}, Hochschild gave a pioneering trial to
the cohomology theory for restricted Lie algebras and as a final
conclusion a 6-term exact sequence was obtained:
\begin{eqnarray*}
0{\to}\H_{\ast}^{1}(L,M)&\to&\H^{1}(L,M)\to
S(L,M^{L})\\
&\to&\H_{\ast}^{2}(L,M)\to \H^{2}(L,M)\to S(L,\H^{1}(L,M)),
\end{eqnarray*}
for $L$  a restricted Lie algebra and $M$  a strongly abelian
restricted Lie algebra with an $L$-operation. Here $\H^{i}$ and
$\H^{i}_{\ast}$ denote the ``ordinary" situation's cohomology groups
and ``restricted" ones respectively,  $S(V,W)$ is the space of
$p$-semilinear maps from $V$ to $W$, and $M^{L}$ is the subset of
invariants, i.e., $M^{L}=\{m\in M|L\cdot m=0\}$. This 6-term exact
sequence established the connection between ordinary cohomology
groups and restricted ones, and was shown to be crucial to get
further information about cohomology theory of restricted Lie
algebras, algebraic groups, infinitesimal groups and discrete groups
\cite{FP1,FP2}. Especially, it can help us to establish Noetherian
property for cohomology algebra $\H^{\ast}(\mathfrak{g},k)$ for a
restricted Lie algebra $\mathfrak{g}$. The current work generalizes
this 6-term exact sequence to restricted Lie superalgebras
\begin{eqnarray*}
0{\to}\H_{\ast}^{1}(\g,M)&{\to}&\H^{1}(\g,M){\to}
S(\g_{\bar{0}},M_{\bar{0}}^{\g})\\
&{\to}&\H_{\ast}^{2}(\g,M){\to} \H^{2}(\g,M){\to}
S(\g_{\bar{0}},\H^{1}(\g,M)).
\end{eqnarray*}
 See Theorem 5.7 for details.  We hope that the result we gotten can be
used as an experimental animal to detect whether the cohomology
algebra of a Lie superalgebra is finitely generated or not.

\subsection{} In Section 2, some necessary notions and results are
collected. In particular, we show that the ordinary (co)homology of
a Lie superalgebra defined by its Koszul complex can be computed
through the Hochschild complex of its enveloping algebra. As we
expect, the extensions of restricted supermodules can be explained
through the first cohomology group $\H^{1}_{\ast}$ . The proof of
this fact is given in Section 3. We give the cohomological
interpretations to the similarity classes and equivalence classes of
extensions of restricted Lie superalgebras in Section 4,5
respectively. And as a result, we get the desired 6-term exact
sequence. The results gotten in the paper are what one would
naturally hope them to be for restricted Lie superalgebras.

Throughout we work with a field $k$ with characteristic $p>2$ as the
ground field. By a superspace we mean a $\mathbb{Z}_{2}$-graded
vector space $V=V_{\overline{0}}\oplus V_{\overline{1}}$, in which
we call elements in $V_{\overline{0}}$ and $V_{\overline{1}}$ even
and odd, respectively. Write $|v|\in \mathbb{Z}_{2}$ for the degree
of $v\in V$, which is implicitly assumed to be
$\mathbb{Z}_{2}$-homogeneous. A linear map $f:\;
V=V_{\overline{0}}\oplus V_{\overline{1}}\to
W=W_{\overline{0}}\oplus W_{\overline{1}}$ is said to be even (resp.
odd) if $f(V_{\bar{i}})\subseteq W_{\bar{i}}$ (resp.
$f(V_{\bar{i}})\subseteq W_{\bar{i+1}}$) for $i=0,1$. Unless
otherwise specified, all vector spaces, algebras, subalgebras,
ideals, modules and submodules etc. are in the super case, and all
linear maps are even. Moreover, for any two $\mathbb{Z}_{2}$-graded
vector spaces $V,W$, we use $\Hom_{k}(V,W)$ to represent the set of
all even linear maps from $V$ to $W$ and $\underline{\Hom}_{k}(V,W)$
to denote that of all  linear maps. A map $f$ from $V$ to $W$ is
$p$-semilinear if $f(\alpha
v_{1}+v_{2})=\alpha^{p}f(v_{1})+f(v_{2})$ for $\alpha\in k$ and
$v_{1},v_{2}\in V$. And, we use the notation $S(V,W)$ to denote the
space of $p$-semilinear maps from $V$ to $W$.

\section{Basic results for the cohomology of (restricted) Lie superalgebras}

The materials in this section are standard generalization from Lie
algebras to Lie superalgebras except Lemma 2.2, where we need give a
generalization of the sign representation of a symmetric group.

\subsection{Basic notions.} The definition of a restricted Lie
superalgebra can be easily formulated (cf. e.g. \cite{WZ}).

\begin{definition}\emph{ A Lie superalgebra $\mathfrak{g}=\mathfrak{g}_{\bar{0}}\oplus
\mathfrak{g}_{\bar{1}}$ is called a \emph{restricted Lie
superalgebra}, if there is a $p$th map $\mathfrak{g}_{\bar{0}}\to
\mathfrak{g}_{\bar{0}}$, denoted as $^{[p]}$, satisfying}

\emph{(a) } \emph{$(cx)^{[p]}=c^{p}x^{[p]}$ for all $c\in k$ and
$x\in \mathfrak{g}_{\bar{0}}$,}

\emph{(b) }\emph{$[x^{[p]},y]=(adx)^{p}(y)$ for all $x\in
\mathfrak{g}_{\bar{0}}$ and $y\in \mathfrak{g}$,}

\emph{(c)
}\emph{$(x+y)^{[p]}=x^{[p]}+y^{[p]}+\sum_{i=1}^{p-1}s_{i}(x,y)$ for
all $x,y\in \mathfrak{g}_{\bar{0}}$ where $is_{i}$ is the
coefficient of $\lambda^{i-1}$ in $(ad(\lambda x+y))^{p-1}(x)$.}
\end{definition}

In short, a restricted Lie superalgebra is a Lie superalgebra whose
even subalgebra is a restricted Lie algebra and the odd part is a
restricted module by the adjoint action of the even subalgebra.  For
a Lie superalgebra $\mathfrak{g}$, $U(\mathfrak{g})$ is denoted to
be its universal enveloping algebra and
$\mathbf{u}(\mathfrak{g})=U(\mathfrak{g})/(x^{p}-x^{[p]}|x\in
\mathfrak{g}_{\bar{0}})$ its restricted enveloping algebra if
moreover $\g$ is restricted.

The notion of the cohomology for a Lie superalgebra
$\mathfrak{g}=\mathfrak{g}_{\bar{0}}\oplus \mathfrak{g}_{\bar{1}}$
was introduced by Fuks (cf. e.g. \cite{F}). By definition, the
(ordinary) space of $n$-dimensional cocycles of $\mathfrak{g}$ with
coefficients in the $\mathfrak{g}$-module $M=M_{\bar{0}}\oplus
M_{\bar{1}}$ is defined to be
$$C^{n}(\mathfrak{g}, M):= \bigoplus_{n_{0}+n_{1}=n}\Hom_{k}(\Lambda^{n_{0}}\mathfrak{g}_{\bar{0}}\otimes
S^{n_{1}}\mathfrak{g}_{\bar{1}}, M).$$ The differential
$\delta_{n}:\;C^{n}(\mathfrak{g}, M)\to C^{n+1}(\mathfrak{g}, M)$ is
defined by the formula\\[2mm]
$\delta_{n-1}(f)(x_{1},\ldots,x_{n_{0}},y_{1},\ldots,y_{n_{1}})$
\begin{eqnarray*}
&=&\sum_{s=1}^{n_{0}}(-1)^{s-1}x_{s}\cdot
f(x_{1},\ldots,\hat{x_{s}},\ldots,x_{n_{0}},y_{1},\ldots,y_{n_{1}})\\
&&+\sum_{t=1}^{n_{1}}(-1)^{n_{0}}y_{t}\cdot
f(x_{1},\ldots,x_{n_{0}},
y_{1},\ldots,\hat{y_{t}},\ldots,y_{n_{1}})\\
&&+\sum_{1\leq s<t\leq
n_{0}}(-1)^{s+t}f([x_{s},x_{t}],x_{1},\ldots,\hat{x_{s}},\ldots,\hat{x_{t}},\ldots,x_{n_{0}},
y_{1},\ldots,y_{n_{1}})\\
&&+\sum_{s=1}^{n_{0}}\sum_{t=1}^{n_{1}}(-1)^{s}f(x_{1},\ldots,\hat{x_{s}},\ldots,x_{n_{0}},[x_{s},y_{t}],
y_{1},\ldots,\hat{y_{t}},\ldots,y_{n_{1}})\\
&&+\sum_{1\leq s<t\leq
n_{1}}-f([y_{s},y_{t}],x_{1},\ldots,x_{n_{0}},
y_{1},\ldots,\hat{y_{s}},\ldots,\hat{y_{t}},\ldots,y_{n_{1}}).
\end{eqnarray*}

For $x_{1},\ldots, x_{n}\in \mathfrak{g}_{\bar{0}}\cup \mathfrak{g}_{\bar{1}}$, one can give a more unified expression for the differential\\[2mm]
$\delta_{n-1}(f)(x_{1},\ldots,x_{n})$
\begin{eqnarray*}
&=&\sum_{s=1}^{n} (-1)^{s-1+|x_{s}|(\sum_{i=1}^{s-1}|x_{i}|)}
x_{s}\cdot
f(x_{1},\ldots,\hat{x_{s}},\ldots,x_{n})\\
&&+\sum_{1\leq s<t\leq
n_{0}}(-1)^{s+t+|x_{s}|(\sum_{i=1}^{s-1}|x_{i}|)+
|x_{t}|(\sum_{i=1}^{t-1}|x_{i}|)+|x_{s}||x_{t}|}\\
&&\cdot
f([x_{s},x_{t}],x_{1},\ldots,\hat{x_{s}},\ldots,\hat{x_{t}},\ldots,x_{n}).
\end{eqnarray*}

It is straightforward to show that $\delta_{n+1}\circ \delta_{n}=0$
and hence, in particular, one can define the cohomologies by setting
$$\H^{n}(\mathfrak{g},M):=\Ker \delta_{n}/\Im \delta_{n-1}.$$
We call it the\emph{ $n$-th cohomology of $\mathfrak{g}$ with
coefficients in $M$}.

Also one can use the usual Hochschild's complex to define the
cohomologies for any augmented algebra. In our case, let
$U(\mathfrak{g})^{+}$ be the ideal in $U(\mathfrak{g})$ generated by
$\mathfrak{g}$. The $n$-cochains are now the \emph{even} $n$-linear
functions on $U(\mathfrak{g})^{+}$ with values in $M$, and the
coboundary operator $\delta$ is defined by the formula
$$\delta_{n-1}(f)(s_{1},\ldots,s_{n})=s_{0}\cdot f(s_{2},\ldots,s_{n})+\sum_{i=1}^{n-1}(-1)^{i}f(s_{1},\ldots,s_{i}s_{i+1},\ldots,s_{n}),$$
for $s_{1},\ldots,s_{n}\in U(\g)^{+}$.

 At the first glance, it is
not so clear whether the cohomologies as defined above in the two
different ways are same or not. For convenience, we call the
cochains defined in the first way and the second way, the \emph{Lie}
type and \emph{associative} type, respectively.

\begin{lemma} There is a canonical isomorphism between the cohomology
groups of Lie type and associative type.
\end{lemma}
\begin{proof} We give an explicit cochain map between two complexes
defined as above. To do it, we need introduce a notation at first.
Let $S_{n}$ be the symmetric group in $n$ letters. For any
$\sigma\in S_{n}$ and $1\leq n_{0}\leq n$, define
$$\sgn(\sigma(n_{0}|n)):=(-1)^{\bar{\sigma(1)}+\cdots +\bar{\sigma(n)}}$$
where $\bar{\sigma(i)}:=\#\{j\in\{1,\ldots,n_{0}\}|j\notin
\{\sigma(1),\ldots,\sigma(i-1)\}, j< \sigma(i)\}$ for $1\leq i\leq
n$. Now, for every cochain $f$ of associative type, define a cochain
$f'$ of Lie type by the formula
\begin{equation}f'(x_{1},\ldots,x_{n})=\sum_{\sigma\in
S_{n}}\sgn(\sigma(n_{0}|n))f(x_{\sigma(1)},\ldots,x_{\sigma(1)})
\end{equation}
where we assume that $x_{1},\ldots,x_{n_{0}}\in
\mathfrak{g}_{\bar{0}}$ while $x_{n_{0}+1},\ldots,x_{n}\in
\mathfrak{g}_{\bar{1}}$. One can verify directly the
$\delta(f')=(\delta(f))'$ and indeed the map $f\mapsto f'$ induces
an isomorphism of the cohomology groups.
\end{proof}

\begin{remark} \emph{The notion $\sgn(\sigma(n_{0}|n))$ generalizes our
usual sign representation $\sgn:\;S_{n}\to \{\pm 1\}$. In fact, we
always have} $$\sgn(\sigma(n|n))=\sgn(\sigma).$$ \emph{Therefore,
the formula (2.1) generalizes the isomorphism given by Hochschild
(\cite{Hoch}, p. 557) for Lie algebra to Lie superalgebra.}
\end{remark}

Note that we can not use Lie type cochains to define the cohomology
groups for a restricted Lie superalgebra directly. Compare with Lie
type cochains, we still can use associative type cochains to define
cohomologies for restricted Lie superalgebras. In this case, we just
need replace $U(\mathfrak{g})^{+}$ by $\mathbf{u}(\mathfrak{g})^{+}$
and $M$ by a restricted $\mathfrak{g}$-module. Similar to the
definition of $U(\mathfrak{g})^{+}$, $\mathbf{u}(\mathfrak{g})^{+}$
is the ideal in $\mathbf{u}(\mathfrak{g})$ generated by
$\mathfrak{g}$. In order to not cause confusion, the restricted
cohomology groups are denoted by
$$\H_{\ast}^{n}(\mathfrak{g},M),\;\;n\in \mathbb{N}.$$ Since the canonical
homomorphism $U(\mathfrak{g})\to \mathbf{u}(\mathfrak{g})$ allows us
to regard any $\mathbf{u}(\mathfrak{g})$-module also as a
$U(\mathfrak{g})$-module, there is a canonical homomorphism
$$\H_{\ast}^{n}(\mathfrak{g},M)\to \H^{n}(\mathfrak{g},M),\;\;n\in \mathbb{N}.$$
An explicit cochain map inducing this homomorphism is given, in the
associative type, by $f\mapsto f^{0}$, where
$f^{0}(x_{1},\ldots,x_{n})=f(x'_{1},\ldots,x'_{n})$ with $x_{i}\in
U(\mathfrak{g})^{+}$ and $x_{i}'$ its canonical image in
$\mathbf{u}(\mathfrak{g})^{+}$.

Denote the complex of the cochains for $U(\mathfrak{g})^{+}$ in the
restricted $\mathfrak{g}$-module $M$ by $C(M)$, and let $C^{0}(M)$
stand for the subcomplex consisting of the cochains of the form
$f^{0}$ with $f$ a cochain for $U(\mathfrak{g})^{+}$ in $M$. Then we
have an exact sequence of complexes $0\to C^{0}(M)\to C(M)\to
C(M)/C^{0}(M)\to 0$ and so we get a long exact sequence
$$\cdots\to \H^{n}_{\ast}(\mathfrak{g},M)\to \H^{n}(\mathfrak{g},M)\to \H^{n}(C(M)/C^{0}(M))\to
\H^{n+1}_{\ast}(\mathfrak{g},M)\to \cdots.$$ Obviously,
$\H^{0}(C(M)/C^{0}(M))=0$ and so there is an injection
\begin{equation}i_{1}:\;\H^{1}_{\ast}(\mathfrak{g},M)\hookrightarrow
\H^{1}(\mathfrak{g},M).\end{equation}

\subsection{Extensions.} Let $\mathfrak{g}$ be a Lie superalgebra
and $K,N$ two $\mathfrak{g}$-modules. An \emph{extension of $K$ by
$N$} is a pair $(E,\phi)$, where $E$ is a $\mathfrak{g}$-module
containing $K$, and $\phi$ is a $\mathfrak{g}$-epimorphism $E \to N$
such that $\Ker \phi=K$. That is, there is an exact sequence of
$\mathfrak{g}$-modules
$$0\to K\to E\stackrel{\phi}{\to} N\to 0.$$
Two such extensions $(E,\phi)$ and $(E',\phi')$ are said to be
\emph{equivalent} if there is a $\g$-isomorphism $\alpha:\;E\to E'$
which leaves the elements of $K$ fixed and satisfies the relation
$\phi'\alpha=\phi$. As usual, denote the equivalence classes of the
extensions of $K$ by $N$ by $\Ext(K,N)$, and there is an ordinary
(i.e., not super) linear space structure over $\Ext(K,N)$. Define
$M$ to be the $k$-space consisting of \emph{all} $k$-linear maps
from $N$ to $K$, that is, using our notion
$$M=\underline{\Hom}_{k}(N,K).$$
Through setting
$M_{\bar{0}}:=\underline{\Hom}_{k}(N_{\bar{0}},K_{\bar{0}})\oplus
\underline{\Hom}_{k}(N_{\bar{1}},K_{\bar{1}}),\;M_{\bar{1}}:=\underline{\Hom}_{k}(N_{\bar{0}},K_{\bar{1}})\oplus
\underline{\Hom}_{k}(N_{\bar{1}},K_{\bar{0}})$, $M=M_{\bar{0}}\oplus
M_{\bar{1}}$ is a superspace. And it is indeed a $\g$-module through
\begin{equation}(x\cdot m)(a):=x\cdot m(a)-(-1)^{|x||m|}m(x\cdot a)\end{equation} for $x\in
\g, m\in M$ and $a\in N$.

\begin{lemma} With notions as above, there is an isomorphism of $k$-spaces
$$\Ext(K,N)\cong \H^{1}(\g, M).$$
\end{lemma}
\begin{proof} For any extension $(E,\phi)$ of $K$ by $N$, let
$\varphi:\; N\to E$ be a linear map such that
$\phi\varphi=\textrm{id}_{N}$. Note that $\varphi$ can be chosen as
an even linear map due to $\phi$ is already an even homomorphism.
From this, we obtain a linear map $f:\; \g \to M$ by setting
$$f(x)(a):=x\cdot \varphi(a)-\varphi(x\cdot a)$$
for $x\in \g, a\in N$. It is straightforward to show that $f$ is a
1-cocycle and above procedure induces a linear map
$F:\;\Ext(K,N)\to\H^{1}(\g, M)$.

Conversely, for any 1-cocycle $f\in \Hom_{k}(\g, M)$. We can attach
it with an extension $(E_{f},\phi_{f})$ of $K$ by $N$. By
definition, as vector space $E_{f}=K\oplus N,\;\phi_{f}(c,d)=d$ and
the $\g$-module is given through
$$x\cdot(c+d):=x\cdot c+ x\cdot d + f(x)(d)$$
for $x\in \g, c\in K$ and $d\in N$. Also, one can show this process
gives a linear map $G:\; \H^{1}(\g, M)\to \Ext(K,N)$. At last, it is
direct to prove that $FG=\textrm{id}_{\H^{1}(\g, M)}$ and
$GF=\textrm{id}_{\Ext(K,N)}$.
\end{proof}

Just like the Lie algebra case, we also hope that we can give an
interpretation to the second cohomology group by using the
extensions of Lie superalgebras. To attack it, let $M$ an abelian
Lie superalgebra, $\g$ an arbitrary Lie superalgebra. An
\emph{extension} \emph{of $M$ by $\g$} is a pair $(E,\phi)$, where
$E$ is a Lie superalgebra containing $M$ as an ideal, and $\phi$ is
a Lie superalgebra epimorphism $E \to \g$ such that $\Ker \phi=K$.
That is, there is an exact sequence of Lie superalgebras
$$0\to M\to E\stackrel{\phi}{\to} \g\to 0.$$ This situation defines
on $M$ the structure of a $\g$-module, with $\g$ operating on $M$
via $E$, in the natural fashion. Similarly, two such extensions
$(E,\phi)$ and $(E',\phi')$ are said to be \emph{equivalent} if
there is a Lie superalgebra isomorphism $\alpha:\;E\to E'$ ($\alpha$
is even by definition) which leaves the elements of $M$ fixed and
satisfies the relation $\phi'\alpha=\phi$. Denote the equivalence
classes of $M$ by $\g$ by $\Ext(M,\g)$ and it is also an ordinary
linear space.

\begin{lemma} With notions as above, there is an isomorphism of
linear spaces
$$\Ext(M,\g)\cong \H^{2}(\g,M).$$
\end{lemma}
\begin{proof} Similar to the proof of Lemma 2.4, we just give the formula for the
correspondence and leave the reader to check the details. For any
extension $(E,\phi)$ of $M$ by $\g$, let $\varphi:\; M\to E$ be a
linear map such that $\phi\varphi=\textrm{id}_{M}$. From this, we
obtain a linear map $f\in\Hom_{k}(\g\otimes \g, M)$ by setting
$$f(x_{1},x_{2}):=[\varphi(x_{1}),\varphi(x_{2})]-\varphi([x_{1},x_{2}])$$
for $x_{1},x_{2}\in \g$. It is straightforward to show that $f$ is a
2-cocycle and above procedure induces a linear map
$F:\;\Ext(M,\g)\to\H^{2}(\g, M)$.

Conversely, for any 2-cocycle $f\in \Hom_{k}(\g\otimes \g, M)$. We
can attach it with an extension $(E_{f},\phi_{f})$ of $M$ by $\g$.
By definition, as vector space $E_{f}=\g\oplus
M,\;\phi_{f}(x_{1},m_{1})=x_{1}$ and the Lie superalgebra structure
is given through
$$[(x_{1},m_{1}),(x_{2},m_{2})]:=([x_{1},x_{2}], x_{1}\cdot m_{1}-(-1)^{|x_{1}||x_{2}|}x_{2}\cdot m_{1}+ f(x_{1},x_{2}))$$
for $x_{1},x_{2}\in \g$ and $m_{1},m_{2}\in M$. Here we implicitly
ask $(x_{i},m_{i})$ to be an homogeneous element and thus we always
have $|x_{i}|=|m_{i}|$. Also, one can show this process gives a
linear map $G:\; \H^{2}(\g, M)\to \Ext(M,\g)$. At last, it is direct
to prove that $FG=\textrm{id}_{\H^{2}(\g, M)}$ and
$GF=\textrm{id}_{\Ext(M,\g)}$.
\end{proof}

\begin{remark} \emph{Both Lemma 2.4 and Lemma 2.5 should be known by experts.
The author just has not found a suitable reference.}
\end{remark}

For latter use and completeness, we collect some identities, which
already appeared in \cite{Hoch}.

\begin{lemma} Let $k\{x,y\}$ be the free algebra generated by
two variables $x,y$. If $D_{w}$ denote the map $z\mapsto
wz-zw=D_{w}(z)$, then we have

\emph{(1)} $\sum_{i=0}^{p-1}x^{i}yx^{p-1-i}=D_{x^{p-1}}(y)$.

\emph{(2)}
$\sum_{i=0}^{l-1}x^{i}D_{x}^{l-1-i}(y)=\sum_{j=0}^{l-1}(-1)^{j}\left
(
\begin{array}{c} l\\j+1
\end{array}\right) x^{l-1-j}yx^{j}.$
\end{lemma}

\begin{proof} (1) It is not hard to see that
$D_{x}(\sum_{i=0}^{p-1}x^{i}yx^{p-1-i}-D_{x^{p-1}}(y))=0$ which
implies $\sum_{i=0}^{p-1}x^{i}yx^{p-1-i}=D_{x^{p-1}}(y)$.

(2) Consider the commutative polynomial ring $k[x_{1},x_{2}]$ at
first. In such a ring, we always have
$$x_{1}^{l}-(x_{1}-x_{2})^{l}=x_{2}\sum_{i=0}^{l-1}x_{1}^{i}(x_{1}-x_{2})^{l-1-i}$$
which implies
$$\sum_{i=0}^{l-1}x_{1}^{i}(x_{1}-x_{2})^{l-1-i}=\sum_{j=0}^{l-1}(-1)^{j}\left
(
\begin{array}{c} l\\j+1
\end{array}\right)x_{1}^{l-1-j}x_{2}^{j}.$$
By specializing this to our case where $x_{1}$ is the left
multiplication by $x$ in $k\{x,y\}$ and $x_{2}$ the right
multiplication by $x$ in $k\{x,y\}$, we get the desired equation.
\end{proof}

\section{Extensions of restricted modules }

In Subsection 2.2, we have considered the extensions of
supermodules. Now let us consider the analogous situation in the
case where $\g$ is a restricted Lie superalgebra, and $K,N$ are
restricted $\g$-modules. Correspondingly, an extension $(E,\phi)$ of
$K$ by $N$ is then called \emph{a restricted extension} if $E$ is a
restricted $\g$-module. Through equation (2.3),
$M=\underline{\Hom}_{k}(N,K)$ is a $\g$-module.

\begin{lemma} $M$ is also a restricted $\g$-module.
\end{lemma}
\begin{proof} To show it, for any $x\in \g_{\bar{0}}$ we define two maps $u_{x},v_{x}\in \Hom_{k}(M,M)$.
By definition, $u_{x}(f)(a):=x\cdot f(a),\;\;v_{x}(f)(a):=f(x\cdot
a)$ for $f\in M$ and $a\in N$. Clearly, $u_{x}v_{x}=v_{x}u_{x}$ and
$x\cdot f=u_{x}(f)-v_{x}(f)$. So,
\begin{eqnarray*}
(x^{p}\cdot
f)(a)&=&(u_{x}-v_{x})^{p}(f)(a)=(u_{x}^{p}-v_{x}^{p})(f)(a)=x^{p}\cdot
f(a)-f(x^{p}\cdot a)\\
&=&x^{[p]}\cdot f(a)-f(x^{[p]}\cdot a)=(x^{[p]}\cdot f)(a)
\end{eqnarray*}
for $x\in \g_{\bar{0}}, f\in M$ and $a\in N$.
\end{proof}

The equivalence classes of restricted extensions of $K$ by $N$ is
denoted by $\Ext_{\ast}(K,N)$. Since any restricted extension can be
regarded as an ordinary extension naturally, there is a natural
linear map $i_{2}:\;\Ext_{\ast}(K,N)\hookrightarrow \Ext(K,N)$. The
main result of this section is the following conclusion.

\begin{proposition} Let $\g$ be a restricted Lie superalgebra, $K,N$ two restricted $\g$-modules, and
 $M=\underline{\Hom}_{k}(N,K)$. Then there is a canonical isomorphism $F|_{\ast}:\; \Ext_{\ast}(K,N) \stackrel{\cong}{\to}
 \H^{1}_{\ast}(\g,M)$ such that the following diagram of canonical maps is commutative

 \begin{figure}[hbt]
\begin{picture}(150,70)(0,0)
\put(0,60){\makebox(0,0){$\Ext_{\ast}(K,N)$}}\put(70,70){\makebox(0,0){$F|_{\ast}$}}
\put(70,10){\makebox(0,0){$F$}} \put(150,60){\makebox(0,0){$
\H^{1}_{\ast}(\g,M)$}} \put(40,60){\vector(1,0){70}}
\put(10,50){\vector(0,-1){40}}\put(20,30){\makebox(0,0){$i_{2}$}}\put(130,30){\makebox(0,0){$i_{1}$}}
\put(140,50){\vector(0,-1){40}}
\put(0,0){\makebox(0,0){$\Ext(K,N)$}} \put(150,0){\makebox(0,0){$
\H^{1}(\g,M)$}} \put(40,0){\vector(1,0){70}}
\end{picture}
\end{figure}

where $F$ is the isomorphism given in Lemma $2.4$ and $i_{1}$ is the
injection described in $(2.2)$.
\end{proposition}
\begin{proof} To attack it, it is enough to show that $Fi_{2}(\Ext_{\ast}(K,N))=\Im
i_{1}$. Actually, we will show that both $Fi_{2}(\Ext_{\ast}(K,N))$
and $\Im i_{1}$ equal to the  subspace $V$ of $\H^{1}(\g,M)$ whose
elements are represented by Lie type 1-cocycles $f$ satisfying
$$x^{p-1}\cdot f(x)=f(x^{[p]})$$ for $ x\in \g_{\bar{0}}.$

Let $(E,\phi)$ be a restricted extension of $K$ by $N$. By regarding
it as a usual extension, we get a 1-cocycle $f\in \Hom_{k}(\g,
\underline{\Hom}_{k}(N,K))$. By definition, $f(x)=x\cdot \varphi$
for a linear section $\varphi$ of $\phi$ and $x\in \g$. Thus
$x^{p-1}\cdot f(x)=x^{p}\cdot \varphi$. Since $E$ and $N$ are
restricted modules, it follows that $x^{p-1}\cdot f(x)=f(x^{[p]})$
for all $x\in \g_{\bar{0}}$. Therefore,
$Fi_{2}(\Ext_{\ast}(K,N))\subseteq V$. Conversely, let $f\in
\Hom_{k}(\g, \underline{\Hom}_{k}(N,K))$ be any Lie type 1-cocycle
satisfying $x^{p-1}\cdot f(x)=f(x^{[p]})$, for $ x\in \g_{\bar{0}}$.
By the construction introduced in the proof of Lemma 2.4, the
corresponding usual extension is denoted by $(E_{f},\phi_{f})$. By
the definition of $E_{f}$, it implies that
$$x^{p}\cdot(c+d)=x^{p}\cdot c+x^{p}\cdot d+\sum_{i=0}^{p-1}x^{i}\cdot f(x)(x^{p-1-i}\cdot d)$$
for $x\in \g_{\bar{0}},\;c\in K$ and $d\in N$. By $K,N$ are
restricted modules and Lemma 2.7 (1), $x^{p}\cdot(c+d)=x^{[p]}\cdot
c+x^{[p]}\cdot d+D_{x^{p-1}}(f(x))=x^{[p]}\cdot c+x^{[p]}\cdot
d+x^{p-1}\cdot f(x)=x^{[p]}\cdot c+x^{[p]}\cdot
d+f(x^{[p]})=x^{[p]}\cdot(c+d)$. Hence, $V\subseteq
Fi_{2}(\Ext_{\ast}(K,N))$ and therefore
$V=Fi_{2}(\Ext_{\ast}(K,N))$.

Let $f$ be Lie type 1-cocycle and assume its cohomology class
belongs to $\Im i_{1}$. Thus there is an associative type 1-cocycle
$g\in \Hom(\mathbf{u}(\g)^{+},\underline{\Hom}_{k}(N,K))$ such that
the cohomology class of $g^{0}$ is same as that of $f$. Thus there
is an element $m\in \underline{\Hom}_{k}(N,K)$ such that
$f(x)=g(x)+x\cdot m$ for $x\in \g$.  Note that $g$ is 1-cocycle,
$g(xy)=x\cdot g(y)$ for any $x,y\in \mathbf{u}(\g)^{+}$. Therefore,
$x^{p-1}\cdot f(x)=x^{p-1}\cdot g(x)+ x^{p}\cdot
m=g(x^{[p]})+x^{[p]}\cdot m=f(x^{[p]})$. So $\Im i_{1}\subseteq V$.
Conversely, let $f\in \Hom_{k}(\g, \underline{\Hom}_{k}(N,K))$ be
any Lie type 1-cocycle satisfying $x^{p-1}\cdot f(x)=f(x^{[p]})$,
for $ x\in \g_{\bar{0}}$. Denote one of its corresponding
associative type 1-cocycles by $g$. Thus there is an element $m\in
\underline{\Hom}_{k}(N,K)$ such that $g(x)=f(x)+x\cdot m$. So
$g(x^{p}-x^{[p]})=g(x^{p})-f(x^{[p]})-x^{[p]}\cdot m=x^{p-1}\cdot
g(x)-f(x^{[p]})-x^{[p]}\cdot m=x^{p-1}\cdot(f(x)+x\cdot
m)-f(x^{[p]})-x^{[p]}\cdot m=0$. So $g$ is indeed defined over
$\mathbf{u}(\g)$. Thus $V\subseteq \Im i_{1}$.
\end{proof}

\section{Extensions of restricted Lie superalgebras: the similarity classes}

The definition of a Lie superalgebra extension has been described in
Subsection 2.2.  We also hope to consider the analogous case where
Lie superalgebras are replaced by restricted ones. The definition of
a restricted Lie superalgebra extension can be given directly. Let
$M$ be an abelian restricted Lie superalgebra and $\g$ just a
restricted Lie superalgebra. A \emph{restricted extension of $M$ by
$\g$} is a pair $(E,\phi)$ where $E$ is a restricted Lie
superalgebra containing $M$ as an ideal, and $\phi$ is a restricted
Lie superalgebra epimorphism $E \to \g$ such that $\Ker \phi=M$.
Similarly, this situation defines on $M$ the structure of a
$\g$-module and it is easy to see this module is restricted.

Obviously, there are two ways to consider the relations between
different restricted extensions: \emph{similarity classes} and
\emph{equivalence classes}. By definition, two restricted extensions
$(E,\phi)$ and $(E',\phi')$ are said to be \emph{similar} if there
is a Lie superalgebra isomorphism $\alpha:\;E\to E'$  which leaves
the elements of $M$ fixed and satisfies the relation
$\phi'\alpha=\phi$. And, they are \emph{equivalent} if moreover
$\alpha$ is a restricted map, that is,
$\alpha(x^{[p]})=\alpha(x)^{[p]}$ for $x\in E_{\bar{0}}$. In this
section, we want to characterize the similarity classes by using
cohomology theory. In subsection 2.2, we have used the notion
$\Ext(M,\g)$ to denote the set of ordinary equivalence classes. To
not cause confusion, we introduce two more notions. The set of
similarity classes and equivalence classes of restricted extensions
of $M$ by $\g$ are denoted by $\Ext_{0}(M,\g)$ and
$\Ext_{\ast}(M,\g)$, respectively. They are abelian groups. Clearly,
we have two natural maps of abelian groups
\begin{equation}i_{3}:\;\Ext_{0}(M,\g)\hookrightarrow
\Ext(M,\g),\;\;\;\;\pi_{1}:\;\Ext_{\ast}(M,\g)\twoheadrightarrow
\Ext_{0}(M,\g)\end{equation} where $i_{3}$ is injective and
$\pi_{1}$ is surjective. As we will see later, both $\Ext_{0}(M,\g)$
and $\Ext_{\ast}(M,\g)$ are ordinary linear spaces whenever $M$ is
\emph{strongly abelian}. And in such case, above two maps are linear
maps automatically.

\begin{definition} A restricted Lie superalgebra $M$ is \emph{strongly
abelian} if it is abelian and $x^{[p]}=0$ for all $x \in
M_{\bar{0}}$.
\end{definition}

\begin{lemma} Every restricted extension of $M$ by $\g$ is similar
to one in which $M$ is strongly abelian.
\end{lemma}
\begin{proof} Note that we can generalize the Proposition 2.1 in
\cite{SF} to the Lie superalgebra directly. That is, for a
restricted Lie superalgebra $(\g, [p])$ and a map $[p]_{1}:\;\g\to
\g$, $[p]_{1}$ is still a $p$-mapping of $\g$ if and only if
$[p]-[p]_{1}$ is a $p$-semilinear map from $\g$ to $C(\g)$, where
$C(\g)$ is the center of $\g$.

Now let $(E,\phi)$ be a restricted extension. By the definition of a
$p$-mapping and $M$ is an abelian ideal of $E$, $M^{[p]}$ is
contained in the center of $E$. So the restriction of $[p]$ to $M$
sends $M$ to $C(E)$ and this map is $p$-semilinear. We evidently
extend this map to a $p$-semilinear map $g$ from $E$ to $C(E)$. Thus
by the result stated in the above paragraph, $[p]_{1}:=[p]-g$ is
another $p$-mapping of $E$. Equipping with this new $p$-mapping, we
get a restricted extension of $M$ by $\g$ which is clearly similar
to the given extension by the identity map, and in which $M$ is
strongly abelian.
\end{proof}

By this lemma, this is no harm to assume that $M$ is strong abelian,
and we indeed do so in the following of this section, when we only
consider similarity classes. By Lemma 2.5, there is an isomorphism
between and $\Ext(M,\g)$ and $\H^{2}(\g, M)$. Since $\Ext_{0}(M,\g)$
is a subset of $\Ext(M,\g)$, there is a subset of $\H^{2}(\g, M)$
which corresponds to $\Ext_{0}(M,\g)$. For convenience, denote this
subset by $\H_{0}^{2}(\g, M)$. So our aim is to characterize
$\H_{0}^{2}(\g, M)$ by using cohomologies.

Now, let $f$ be a 2-cocycle. Recall from the proof of Lemma 2.5 the
corresponding extension is $(E_{f},\phi_{f})$. For $x\in
\g_{\bar{0}}, x_{1}\in \g$, direct computations show that
$$D_{(x,0)^{p}}(x_{1},0)=([x^{[p]},x_{1}], \sum_{i=0}^{p-1}x^{i}\cdot f(x,D_{x^{p-1-i}}(x_{1}))).$$
For short, define $k_{x}(x_{1}):=\sum_{i=0}^{p-1}x^{i}\cdot
f(x,D_{x^{p-1-i}}(x_{1}))$ and $f_{x_{1}}(x):=f(x,x_{1})$.

\begin{lemma} The map $k_{x}+f_{x^{[p]}}$ for any $x\in \g_{\bar{0}}$  is a 1-cocycle from $\g$
to $M$ and it only depends on the cohomology class of $f$.
\end{lemma}
\begin{proof} By $x\in \g_{\bar{0}}$, $D_{(x,0)^{p}}$ is an ordinary
derivation, that is
$$D_{(x,0)^{p}}[(x_{1},0),(x_{2},0)]=[D_{(x,0)^{p}}(x_{1},0),
(x_{2},0)]+[(x_{1},0), D_{(x,0)^{p}}(x_{2},0)]$$ for $x_{1},x_{2}\in
\g$. From this, we have
\begin{eqnarray*}
k_{x}([x_{1},x_{2}])+x^{[p]}\cdot f(x_{1},x_{2})&=&x_{1}\cdot
k_{x}(x_{2})-(-1)^{|x_{1}||x_{2}|}x_{2}\cdot k_{x}(x_{1})\\
&+& f(x_{1},[x^{[p]},x_{2}])+f([x^{[p]},x_{1}],x_{2}).
\end{eqnarray*}
Since $f$ is a 2-cocycle, $x^{[p]}\cdot
f(x_{1},x_{2})-f(x_{1},[x^{[p]},x_{2}])-f([x^{[p]},x_{1}],x_{2})=-x_{1}\cdot
f(x_{2},x^{[p]})+(-1)^{|x_{1}||x_{2}|}x_{2}\cdot
f(x_{1},x^{[p]})+f([x_{1},x_{2}],x^{[p]})= -x_{1}\cdot
f_{x^{[p]}}(x_{2})+(-1)^{|x_{1}||x_{2}|}x_{2}\cdot
f_{x^{[p]}}(x_{1})+f_{x^{[p]}}([x_{1},x_{2}])$. From this, it is
easy to see that $k_{x}+f_{x^{[p]}}$ is indeed a 1-cocycle. To show
the second claim, we need show $k_{x}+f_{x^{[p]}}$ is a coboundary
whenever $f$ is so. Now assume that $f=\delta g$ for some $g\in
\Hom_{k}(\g,M)$. Then
\begin{eqnarray*} k_{x}(x_{1})&=&\sum_{i=0}^{p-1}x^{i}\cdot
f(x,D_{x^{p-1-i}}(x_{1}))\\
&=&\sum_{i=0}^{p-1}x^{i}\cdot(x\cdot
g(D_{x^{p-1-i}}(x))-D_{x^{p-1-i}}(x_{1})\cdot
g(x)-g(D_{x^{p-i}}(x_{1})))\\
&=& x^{[p]}\cdot g(x_{1})-g([x^{[p]},x_{1}])-
\sum_{i=0}^{p-1}x^{i}D_{x^{p-1-i}}(x_{1})\cdot g(x)\\
&=&x^{[p]}\cdot
g(x_{1})-g([x^{[p]},x_{1}])-\sum_{i=0}^{p-1}(-1)^{i}\left (
\begin{array}{c} p\\i+1
\end{array}\right)x^{p-1-i}x_{1}x^{i}\cdot g(x)\\
&=&x^{[p]}\cdot g(x_{1})-g([x^{[p]},x_{1}])-x_{1}x^{p-1}\cdot g(x),
\end{eqnarray*}
where Lemma 2.7 (2) is used. Therefore,
$k_{x}(x_{1})+f_{x^{[p]}}(x_{1})=x^{[p]}\cdot
g(x_{1})-g([x^{[p]},x_{1}])-x_{1}x^{p-1}\cdot
g(x)+x_{1}\cdot(g(x^{[p]}))-x^{[p]}\cdot
g(x_{1})-g([x_{1},x^{[p]}])=x_{1}\cdot(g(x^{[p]})-x^{p-1}\cdot
g(x))$. So $k_{x}+f_{x^{[p]}}$ is a coboundary too.
\end{proof}

By this lemma, for any $x\in \g_{\bar{0}}$, we get a  map
$\Phi_{x}:\;\H^{2}(\g, M)\to \H^{1}(\g,M)$ which is induced by
$f\mapsto k_{x}+f_{x^{[p]}}$. The cohomology class of $f$ is denoted
by $c(f)$ and we want to give a more controllable representative to
$\Phi_{x}(c(f))$. For this, let $g$ be an associative type 2-cocycle
whose cohomology class is $c(f)$. So we can assume
$$f(x_{1},x_{2})=g(x_{1},x_{2})-(-1)^{|x_{1}||x_{2}|}g(x_{2},x_{1})$$
for $x_{1},x_{2}\in \g$. Using this associative type 2-cocycle and
noting $x\in \g_{\bar{0}}$,
\begin{eqnarray*} k_{x}(x_{1})&=&\sum_{i=0}^{p-1}x^{i}\cdot
f(x,D_{x^{p-1-i}}(x_{1}))\\
&=&\sum_{i=0}^{p-1}x^{i}\cdot(g(x,D_{x^{p-1-i}}(x_{1}))-g(D_{x^{p-1-i}}(x_{1}),x))\\
&=&\sum_{i=0}^{p-1}g(x^{i+1},D_{x^{p-1-i}}(x_{1}))-\sum_{i=0}^{p-1}g(x^{i},
xD_{x^{p-1-i}}(x_{1}))\\
&&-\sum_{i=0}^{p-1}g(x^{i}D_{x^{p-1-i}}(x_{1}),x)+\sum_{i=0}^{p-1}g(x^{i},
D_{x^{p-1-i}}(x_{1})x)\\
&=&\sum_{i=0}^{p-1}g(x^{i+1},D_{x^{p-1-i}}(x_{1}))-\sum_{i=0}^{p-1}g(x^{i},D_{x^{p-i}}(x_{1}))\\
&&-\sum_{i=0}^{p-1}g(x^{i}D_{x^{p-1-i}}(x_{1}),x)\\
&=&g(x^{p},x_{1})-\sum_{i=0}^{p-1}g(x^{i}D_{x^{p-1-i}}(x_{1}),x)\\
&=&g(x^{p},x_{1})-\sum_{i=0}^{p-1}g((-1)^{i}\left (
\begin{array}{c} p\\i+1
\end{array}\right)x^{p-1-i}x_{1}x^{i},x)\\
&=&g(x^{p},x_{1})-g(x_{1}x^{p-1},x)\\
&=&g(x^{p},x_{1})-g(x_{1},x^{p})-x_{1}\cdot g(x^{p-1},x).
\end{eqnarray*}
Therefore,
$$k_{x}(x_{1})+f_{x^{[p]}}(x_{1})=g(x^{p}-x^{[p]},x_{1})-g(x_{1},x^{p}-x^{[p]})-x_{1}\cdot g(x^{p-1},x).$$
Hence, $\Phi_{x}(c(f))$ has a representative 1-cocycle $g_{x}':\;
x_{1}\mapsto g(x^{p}-x^{[p]},x_{1})-g(x_{1},x^{p}-x^{[p]})$. By this
expression, we get a $p$-semilinear map $g':\;\g_{\bar{0}}\to
\H^{1}(\g,M)$ induced by $x \mapsto g'_{x}$. Thus we get a linear
map \begin{equation}\Phi:\;\H^{2}(\g,M)\to
S(\g_{\bar{0}},\H^{1}(\g,M)),\;\;\;\;g\mapsto g'.\end{equation}

 The
characterization of $\H_{0}^{2}(\g,M)$ is described as follows.

\begin{proposition} $\H_{0}^{2}(\g,M)=\Ker \Phi$.
\end{proposition}
\begin{proof} Let $f$ be a 2-cocycle and assume
$c(f)\in \H_{0}^{2}(\g,M)$. Thus the corresponding extension
$(E_{f},\phi_{f})$ is a restricted extension. So for any $x\in
\g_{\bar{0}}$, there is an element $\rho(x)\in M$ such that
$(x,0)^{[p]}=(x^{[p]},\rho(x))$. From
$D_{(x,0)^{p}}=D_{(x,0)^{[p]}}$, we must have
$k_{x}(x_{1})=f(x^{[p]},x_{1})-x_{1}\cdot \rho(x)$, i.e.,
$k_{x}+f_{x^{[p]}}=\delta(-\rho(x))$, so that $\Phi_{x}(c(f))=0$.
Therefore $c(f)\in \Ker \Phi$ and thus $\H_{0}^{2}(\g,M)\subseteq
\Ker \Phi$.

Conversely, assume $c(f)\in \Ker \Phi$. Recall we use $g$ to denote
the associative type 2-cocycle of $f$ such that
$f(x_{1},x_{2})=g(x_{1},x_{2})-(-1)^{|x_{1}||x_{2}|}g(x_{2},x_{1})$.
Thus we get a $p$-semilinear map
$$\sigma:\; \g_{\bar{0}}\to M_{\bar{0}}$$ such that
$g_{x}'(x_{1})=x_{1}\cdot \sigma(x)$. Now we define a $p$-mapping on
$E_{f}$ through \begin{equation} (x,m)^{[p]}:=(x^{[p]},x^{p-1}\cdot
m+g(x^{p-1},x)-\sigma(x)) \end{equation}
 for $x\in \g_{\bar{0}},
m\in M_{\bar{0}}$. Of course, one can show directly (4.3) indeed
gives a $p$-mapping on $E_{f}$ and thus $E_{f}$ is a restricted Lie
superalgebra. Also, one can copy the same computations used in
restricted Lie algebra case (see p. 568-569 in \cite{Hoch}) to show
(4.3) satisfy all conditions of a $p$-mapping. In one word, $c(f)\in
\H_{0}^{2}(\g,M)$ and thus $ \Ker \Phi \subseteq \H_{0}^{2}(\g,M)$.
\end{proof}

\begin{remark} \emph{By this proposition, if $M$ is strongly abelian, we know that $\Ext_{0}(M,\g)$ is
also an ordinary vector space and the canonical map $i_{3}$ given in
(4.1) is a linear map.}
\end{remark}

\section{Extensions of restricted Lie superalgebras: the equivalence classes}

Further, we consider the restricted equivalence classes
$\Ext_{\ast}(M.\g)$ in this section. And, as the final conclusion,
the Hochschild's 6-term exact sequence will be given. As the
beginning, a decomposition of a similarity class into equivalence
classes will be given.

\subsection{Decomposition of similarity classes.} Let $(E,\phi)$ be a restricted extension of $M$ by $\g$
and denote its similarity class by $c$. We want to decompose $c$
into a set $S_{c}$ of equivalence classes.
 For any other representative object $(E',\phi')$ of $c$, there is
a similarity isomorphism $\gamma:\;(E,\phi)\to (E',\phi')$.

\begin{lemma} For any $e\in E_{\bar{0}}$,
$\gamma(e^{[p]})-(\gamma(e))^{[p]}$  depends only on $\phi(e)$.
\end{lemma}
\begin{proof} To attack it, it is enough to show that $\gamma(e^{[p]})-(\gamma(e))^{[p]}=
\gamma((e+m)^{[p]})-(\gamma(e+m))^{[p]}$ for any $m\in M_{\bar{0}}$.
This is just a direct computation.
$$\gamma((e+m)^{[p]})=\gamma(e^{[p]}+m^{[p]}+e^{p-1}\cdot m)=\gamma(e^{[p]})+m^{[p]}+e^{p-1}\cdot m,$$
$$(\gamma(e+m))^{[p]}=(\gamma(e)+m)^{[p]}=\gamma(e)^{[p]}+m^{[p]}+\gamma(e)^{p-1}\cdot m.$$
Note that we always have $e^{p-1}\cdot m=\gamma(e)^{p-1}\cdot m$, we
get $\gamma(e^{[p]})-(\gamma(e))^{[p]}=
\gamma((e+m)^{[p]})-(\gamma(e+m))^{[p]}$.
\end{proof}

By this lemma, for any $e\in E_{\bar{0}}$, one can denote the
difference $\gamma(e^{[p]})-(\gamma(e))^{[p]}$ by $g(\phi(e))$ and
hence we get a map
$$g:\; \g_{\bar{0}}\to M_{\bar{0}}.$$

\begin{lemma} $g$ is a $p$-semilinear map from $\g_{\bar{0}}$ to $
M_{\bar{0}}^{\g}$ where $M_{\bar{0}}^{\g}:=\{m\in M_{\bar{0}}|x\cdot
m=0,\;x\in \g \}$.
\end{lemma}
\begin{proof} By $\gamma$ is a Lie superalgebra map,
$$\gamma(\sum_{i=0}^{p-1}s_{i}(e_{1},e_{2}))=\sum_{i=0}^{p-1}s_{i}(\gamma(e_{1}),\gamma(e_{2}))$$
for  $e_{1},e_{2}\in E_{\bar{0}}$ (See condition (c) of a
$p$-mapping for the definition of $s_{i}(x,y)$). This indeed implies
that $g$ is a $p$-semilinear map.

Furthermore, for any $x=\phi(z)\in \g$, we have
\begin{eqnarray*}
x\cdot g(\phi(e))&=&\gamma(z)\cdot
(\gamma(e^{[p]})-(\gamma(e))^{[p]})\\
&=&\gamma([z,e^{[p]}])+D_{\gamma(e)^{p}}(\gamma(z))\\
&=&\gamma([z,e^{[p]}])+\gamma (D_{e^{p}}(z))\\
&=&0.
\end{eqnarray*}
\end{proof}

As stated in the first paragraph of the proof of Lemma 4.2, one can
give a new $p$-mapping for $E$ by setting $e^{(p)}:=\;
e^{[p]}-g(\phi(e))$ and thus we get a new restricted Lie
superalgebra which gives an extension of $M$ by $\g$. Denote it by
$(E_{g},\phi)$. Now
$\gamma(e^{(p)})=\gamma(e^{[p]})-g(\phi(e))=\gamma(e)^{[p]}$ and so
$(E_{g},\phi)$ is equivalent to $(E',\phi')$. Conversely, for any
$g\in S(\g_{\bar{0}},M_{\bar{0}}^{\g})$, $(E_{g},\phi)$ is similar
to $(E,\phi)$. Moreover, we get an action
$$g^{\ast}:\;\Ext_{\ast}(M,\g)\to \Ext_{\ast}(M,\g),\;\;\;\;(E,\phi)\mapsto(E_{g},\phi).$$
Such discussions give the following basic fact.

\begin{lemma}  Through the map $g\mapsto g^{\ast}$,
$S(\g_{\bar{0}},M_{\bar{0}}^{\g})$ operates transitively on each
$S_{c}$.
\end{lemma}

We hope to determine the kernel of the map $g\mapsto g^{\ast}$. For
this, we need give a cohomology explanation to the automorphisms of
ordinary extensions. Let $(F,\psi)$ be an ordinary extension of $M$
by $\g$, where "ordinary" means the extension need not to be
restricted. An automorphism of $(F,\psi)$ is an  isomorphism of Lie
superalgebras $\alpha:\; F\to F$ which leaves the elements of $M$
fixed and satisfies that relation $\psi\alpha=\psi$. Since
$\alpha(e)-e=\alpha(e+m)-(e+m)$ for any $m\in M$, $\alpha(e)-e$ only
depends on $\psi(e)$ and denote it by $h(\psi(e))$. From this we get
an even linear map $h:\; \g\to M,\;\psi(e)\mapsto h(\psi(e))$ for
$e\in F$.

\begin{lemma} The $h\in \Hom(\g,M)$ defined above is a 1-cocycle.
\end{lemma}
\begin{proof} For any $e_{1},e_{2}\in F$, we always have
$\alpha([e_{1},e_{2}])=[\alpha(e_{1}),\alpha(e_{2})]$. From this, we
have
$h([\psi(e_{1}),\psi(e_{2})])=[e_{1},h(\psi(e_{2}))]-(-1)^{|e_{1}||e_{2}|}[e_{2},h(\psi(e_{1}))]
=[\psi(e_{1}),h(\psi(e_{2}))]-(-1)^{|e_{1}||e_{2}|}[\psi(e_{2}),h(\psi(e_{1}))]$.
This implies that $h$ is a 1-cocycle.
\end{proof}

Conversely, for any 1-cocycle $h$ one can get an automorphism of
$(F,\psi)$ by setting $\alpha:\; F\to F,\; e\mapsto e+h(\psi(e))$.

Now we go back to determine the kernel of the map $g\mapsto
g^{\ast}$.  Assume $g\in S(\g_{\bar{0}},M_{\bar{0}}^{\g})$ is one
lying in the kernel. So $(E_{g},\phi)$ is restricted equivalent to
$(E,\phi)$ for any restricted extension $(E,\phi)$. Let $\gamma: \;
E\to E_{g}$ be the isomorphism. By forgetting the restricted
structure, $\gamma$ gives an automorphism of $(E,\phi)$. Owing to
Lemma 5.4, $\gamma(e)=e+h(\phi(e))$ for some 1-cocycle $h$. Then by
$\gamma(e^{[p]})=\gamma(e)^{[p]}$ for $e\in E_{\bar{0}}$,
$g(\phi(e))=e^{p-1}\cdot
h(\phi(e))+h(\phi(e))^{[p]}-h(\phi(e)^{[p]})$. That is, for any
$x\in \g_{\bar{0}}$, we have
$$g(x)=x^{p-1}\cdot
h(x)+h(x)^{[p]}-h(x^{[p]}).$$ Define $h'(x):=x^{p-1}\cdot
h(x)+h(x)^{[p]}-h(x^{[p]})$ and it is not hard to see that $h'\in
S(\g_{\bar{0}},M_{\bar{0}}^{\g})$. So we get a liner map
\begin{equation} \Psi:\;Z^{1}(\g,M)\to
S(\g_{\bar{0}},M_{\bar{0}}^{\g}),\;\;h\mapsto h',
\end{equation}
where as usual $Z^{1}(\g,M)$ is the space of 1-cocycles for $\g$ in
$M$. Now, we know that $\Im \Psi$ is just the kernel of the map
$g\mapsto g^{\ast}$. Note that the procedure to determine the kernel
does not depends on the choice of equivalence class of $(E,\phi)$.
So we can choose it to be the trivial extension $s_{0}$, that is,
the 0-element in $\Ext_{\ast}(M,\g)$. Define
$$G_{s_{0}}:\;S(\g_{\bar{0}},M_{\bar{0}}^{\g})\to \Ext_{\ast}(M,\g),\;\;g\mapsto g^{\ast}(s_{0}).$$
So $\Im \Psi=\Ker G_{s_{0}}$. Thus, the following 4-term exact
sequence is gotten.

\begin{proposition} Let $M$ be an abelian restricted Lie
superalgebra on which the restricted Lie superalgebra $\g$ operates.
Then we have the following $4$-term exact sequence of abelian groups
$$Z^{1}(\g,M)\stackrel{\Psi}{\to}
S(\g_{\bar{0}},M_{\bar{0}}^{\g})\stackrel{G_{s_{0}}}{\to}\Ext_{\ast}(M,\g)\stackrel{\pi_{1}}{\to}\Ext_{0}(M,\g)\to
0.$$
\end{proposition}

\subsection{Equivalence classes.}
In this subsection, we always assume that $M$ is strongly abelian.
In such case, the connection between $\Ext_{\ast}(M,\g)$ and
$\H_{\ast}^{2}(\g,M)$ is nice.

\begin{proposition} Let $M$ be a strongly abelian restricted Lie
superalgebra on which the restricted Lie superalgebra $\g$ operates.
Then there is a canonical isomorphism $F|_{\ast}:\;
\Ext_{\ast}(M,\g) \stackrel{\cong}{\to}
 \H^{2}_{\ast}(\g,M)$ such that the following diagram of canonical maps is commutative
 \begin{figure}[hbt]
\begin{picture}(150,70)(0,0)
\put(0,60){\makebox(0,0){$\Ext_{\ast}(M,\g)$}}\put(70,70){\makebox(0,0){$F|_{\ast}$}}
\put(70,10){\makebox(0,0){$F|_{0}$}} \put(150,60){\makebox(0,0){$
\H^{2}_{\ast}(\g,M)$}} \put(40,60){\vector(1,0){70}}
\put(10,50){\vector(0,-1){40}}\put(20,30){\makebox(0,0){$\pi_{1}$}}
\put(140,50){\vector(0,-1){40}}
\put(0,0){\makebox(0,0){$\Ext_{0}(M,\g)$}}
\put(150,0){\makebox(0,0){$ \H_{0}^{2}(\g,M)$}}
\put(40,0){\vector(1,0){70}}
\end{picture}
\end{figure}

where $F|_{0}$ is the restriction of the canonical isomorphism
$F:\;\Ext(M,\g) \stackrel{\cong}{\to}
 \H^{2}(\g,M)$ to $\Ext_{0}(M,\g)$.
\end{proposition}

\begin{proof} The idea is similar to that used in the proof
of Lemma 2.5. Given a restricted associative type 2-cocycle $g$ of
$\mathbf{u}(\g)$ with values in $M$. Construct the corresponding
restricted extension $(E_{g}, \phi_{g})$ as follows: As a space,
$E_{g}=\g\oplus M,\;\phi_{g}(x_{1},m_{1})=x_{1}$ and the restricted
Lie superalgebra structure is given through
\begin{eqnarray*} [(x_{1},m_{1}),(x_{2},m_{2})]&:=&([x_{1},x_{2}], x_{1}\cdot
m_{1}-(-1)^{|x_{1}||x_{2}|}x_{2}\cdot m_{1}\\
 &+& g(x_{1},x_{2})-(-1)^{|x_{1}||x_{2}|}g(x_{2},x_{1})),\end{eqnarray*}
 $$(x,m)^{[p]}:=(x^{[p]},x^{p-1}\cdot m+g(x^{p-1},x))$$
for $x_{1},x_{2}\in \g, x\in \g_{\bar{0}}$ and $m_{1},m_{2}\in M,
m\in M_{\bar{0}}$. It is straightforward to show $(E_{g}, \phi_{g})$
is indeed a restricted extension of $M$ by $\g$. Next, we need to
show that the equivalence class of $(E_{g}, \phi_{g})$ depends only
on the cohomology class in $\H^{2}_{\ast}(\g,M)$ of $g$. Let $h$ be
a 1-cochain and we will show that $(E_{g+\delta(h)},
\phi_{g+\delta(h)})$ is equivalent to $(E_{g}, \phi_{g})$. In fact,
define
$$\alpha:\;E_{g+\delta(h)}\to E_{g},\;\;(x,m)\mapsto (x, m+h(x))$$
and direct computations show that $\alpha$ is an equivalence
isomorphism. Thus the map $g\mapsto (E_{g}, \phi_{g})$ induces a
linear homomorphism $$G|_{\ast}:\; \H^{2}_{\ast}(\g,M)\to
\Ext_{\ast}(M,\g).$$

 Conversely, let $(E, \phi)$ be a restricted extension of $M$ by
 $\g$ and $\phi$ can be extended uniquely to a
 homomorphism $\phi'$ from $\mathbf{u}(E)$ to $\mathbf{u}(\g)$. Clearly, $\Ker(\phi')=\mathbf{u}(E)M$.
  Note that we can not apply the same method used in the proof
 of Lemma 2.5 directly since it only gives us a 2-cocycle of Lie
 type. And in the restricted case, we only have associative type
 cochains. To overcome this difficulty, a linear map is needed. It
 is known that $M$ is a  $\mathbf{u}(E)$-module. We claim the the
 identity map of $M$ to $M$ can be extended in one and only way to a
 $\mathbf{u}(E)$-homomorphism from $\mathbf{u}(E)M$  to $M$, each
 being regarded as a $\mathbf{u}(E)$-module in the natural fashion.
 Actually, the map $$\gamma:\;\mathbf{u}(E)M\to M,\;\;\sum_{i}u_{i}m_{i}\mapsto \sum_{i}\phi'(u_{i})\cdot m_{i}$$
 is the desired one. Now, we go back to construct an associative
 type 2-cocycle. To attack it, let $\psi$ be a linear map from $\g$
 to $E$ which is inverse to $\phi$. We can extend $\psi$ to a linear
 map $\psi':\; \mathbf{u}(\g)\to \mathbf{u}(E)$ such that
 $\phi'\psi'=\textrm{id}_{\mathbf{u}(\g)}$. Now define
 $$g:\;\mathbf{u}(\g)^{+}\otimes \mathbf{u}(\g)^{+}\to M,\;\;(x,y)\mapsto \gamma(\psi'(x)\psi'(y)-\psi'(xy)).$$
 One can check directly that $g$ is an associative type 2-cocycle of
 $\mathbf{u}(\g)$ in $M$ and this induces a linear map $$F|_{\ast}:\; \Ext_{\ast}(M,\g) {\to}
 \H^{2}_{\ast}(\g,M).$$  Once we can show that $(E,\phi)$ is
 equivalent to $(E_{g},\phi_{g})$, then
 $G|_{\ast}\circ F|_{\ast}=\textrm{id}_{\Ext_{\ast}(M,\g)}$.
 Actually, define $$\alpha:\;E_{g}\to E,\;\;(x,m)\mapsto
 \psi(x)+m.$$ Through direct computations, we get
 \begin{eqnarray*}
\alpha([(x_{1},m_{1}),(x_{2},m_{2})])&=&\alpha([x_{1},x_{2}],x_{1}\cdot
m_{2}-(-1)^{|x_{1}||x_{2}|}x_{2}\cdot
m_{1}\\&&+g(x_{1},x_{2})-(-1)^{|x_{1}||x_{2}|}g(x_{2},x_{1}))\\
&=&\psi([x_{1},x_{2}])+x_{1}\cdot
m_{2}-(-1)^{|x_{1}||x_{2}|}x_{2}\cdot m_{1}\\&&
+\gamma([\psi(x_{1}),\psi(x_{2})]-\psi([x_{1},x_{2}]))\\
&=& x_{1}\cdot m_{2}-(-1)^{|x_{1}||x_{2}|}x_{2}\cdot m_{1}
+[\psi(x_{1}),\psi(x_{2})]\\
&=&[\alpha(x_{1},m_{1}),\alpha(x_{2},m_{2})].
\end{eqnarray*}
Furthermore, $$\alpha(x,m)^{[p]}=\psi(x)^{[p]}+x^{p-1}\cdot m,$$
$$\alpha((x,m)^{[p]})=\psi(x^{[p]})+x^{p-1}\cdot m+g(x^{p-1},x).$$
And,
$g(x^{p-1},x)=\gamma(\psi'(x^{p-1})\psi(x)-\psi'(x^{p}))=\gamma(\psi(x)^{p}-\psi(x^{[p]}))=\psi(x)^{[p]}-\psi(x^{[p]})$.
So $\alpha(x,m)^{[p]}=\alpha((x,m)^{[p]})$ and hence $\alpha$ is an
equivalence isomorphism.

At last, let us show that $G|_{\ast}$ is a monomorphism. Suppose
that $(E_{g},\phi_{g})$ is a trivial extension. By definition, there
is a homomorphism $\psi:\; \g\to E_{g}$ of restricted Lie
superalgebras such that $\phi\psi=\textrm{id}_{\g}$. Write
$\psi(x_{1})=(x_{1}, -h(x_{1}))$ for $x_{1}\in \g$. By $\psi$ is a
homomorphism of restricted Lie superalgebras, we get
\begin{equation}
g(x_{1},x_{2})-(-1)^{|x_{1}||x_{2}|}g(x_{2},x_{1})=x_{1}\cdot
h(x_{2})-(-1)^{|x_{1}||x_{2}|}x_{2}\cdot h(x_{1})-h([x_{1},x_{2}]),
\end{equation}
\begin{equation}
g(x^{p-1},x)=x^{p-1}\cdot h(x)-h(x^{[p]})
\end{equation}
for $x_{1},x_{2}\in \g$ and $x\in \g_{\bar{0}}$. The equation (5.2)
implies that Lie type cocycle corresponding $g$ is a coboundary. So
$g^{0}$ is also a coboundary (see subsection 2.2 for the definition
of $g^{0}$). That is, there is a 1-cochain $\omega$ for $U(\g)^{+}$
in $M$ such that $g^{0}(u,v)=u\cdot \omega(v)-\omega(uv)$ for
$u,v\in U(\g)^{+}$. So it is not hard to see $\omega|_{\g}-h$ is a
Lie type 1-cocycle and therefore coincides $\varpi|_{\g}$ for an
associative type 1-cocycle $\varpi$ for $U(\g)^{+}$ in $M$.
Replacing $\omega$ by $\omega-\varpi$, one can assume the
$\omega|_{\g}=h$. So (5.3) implies that $g(x^{p-1},x)=x^{p-1}\cdot
\omega(x)-\omega(x^{[p]})$. At the same time,
$g(x^{p-1},x)=g^{0}(x^{p-1},x)=x^{p-1}\cdot
\omega(x)-\omega(x^{p})$. Thus $\omega(x^{[p]})=\omega(x^{p})$ and
so $\omega=f^{0}$ with $f$ a 1-cochain for $\mathbf{u}(\g)^{+}$ in
$M$. Therefore, $g=\delta f$ and we get the desired conclusion.
\end{proof}

\subsection{The 6-term exact sequence.} Now the Hochschid's 6-term
exact sequence is a direct consequence of conclusions we built.

\begin{theorem} Let $M$ be a strongly abelian restricted Lie
superalgebra on which the restricted Lie superalgebra $\g$ operates.
Then we have the following 6-term exact sequence of ordinary linear
spaces
\begin{eqnarray*}
0{\to}\H_{\ast}^{1}(\g,M)&\stackrel{i_{1}}{\To}&\H^{1}(\g,M)\stackrel{\bar{\Psi}}{\To}
S(\g_{\bar{0}},M_{\bar{0}}^{\g})\\
&\stackrel{F|_{\ast}\circ
G_{s_{0}}}{\To}&\H_{\ast}^{2}(\g,M)\stackrel{F|_{0}\circ
\pi_{1}\circ F|^{-1}_{\ast}}{\To} \H^{2}(\g,M)\stackrel{\Phi}{\To}
S(\g_{\bar{0}},\H^{1}(\g,M)).
\end{eqnarray*}
\end{theorem}
\begin{proof} Recall the definition of $\Psi:\;Z^{1}(\g,M)\to
S(\g_{\bar{0}},M_{\bar{0}}^{\g})$, one can find $\Psi$ maps each
1-coboundary to 0 whenever $M$ is strongly abelian. So $\Psi$
induces a linear map $\bar{\Psi}:\;\H^{1}(\g,M)\to
S(\g_{\bar{0}},M_{\bar{0}}^{\g})$ naturally. By combining the
description of $\H_{\ast}^{1}(\g,M)$ given in the proof of
Proposition 3.2, we have the following exact sequence
$$0{\to}\H_{\ast}^{1}(\g,M)\stackrel{i_{1}}{\To}\H^{1}(\g,M)\stackrel{\bar{\Psi}}{\To}
S(\g_{\bar{0}},M_{\bar{0}}^{\g}).$$ Through using Proposition 5.5,
we get the exact sequence
$$0{\to}\H_{\ast}^{1}(\g,M)\to\H^{1}(\g,M)\to
S(\g_{\bar{0}},M_{\bar{0}}^{\g})\to \Ext_{\ast}(M,\g)\to
\Ext_{0}(M,\g)\to 0.$$ By combining the descriptions of
$\Ext_{0}(M,\g)$ and $\Ext_{\ast}(M,\g)$ given in Propositions 4.4
and 5.6 respectively, the desired 6-term exact sequence is followed.
\end{proof}

\begin{remark} \emph{(1) In page 575 in \cite{Hoch}, Hochschild gave the 6-term exact
sequence in the following way}
\begin{eqnarray*}
0{\to}\H_{\ast}^{1}(L,M)&\to&\H^{1}(L,M)\to
S(L,M^{L})\\
&\to&\H_{\ast}^{2}(L,M)\to \H^{2}(L,M)\to S(L,\H^{1}(L,M)),
\end{eqnarray*}
\emph{where $L$ is a restricted Lie algebra and $M$ is a strongly
abelian restricted Lie algebra with an $L$-operation. Clearly, if we
take $\g$ in Theorem 5.7 to be a restricted Lie algebra, then we
recover the original Hochschild's 6-term exact sequence very well.}

\emph{(2) As we have seen, the proof of main result depends on the
interpretations of cohomology groups by using various kinds of
extensions. It is hopeful that one can get the same result or find
applications by filtering the associative cochain complex for
$U(\g)$ in $M$ relative the ideal, which generated by
$x^{p}-x^{[p]}$ for $x\in \g_{\bar{0}}$, and considering the
corresponding spectral sequence. This procedure should relate the
works in \cite{FP1,FP2} to super case.}
\end{remark}

\section*{Acknowledgements} The author is supported by
Japan Society for the Promotion of Science under the item ``JSPS
Postdoctoral Fellowship for Foreign Researchers" and Grant-in-Aid
for Foreign JSPS Fellow. I would gratefully acknowledge JSPS. The
work is also supported by Natural Science Foundation (No. 10801069).
I would like thank Professor A. Masuoka for stimulating discussions
and his encouragements.


\begin{thebibliography}{99}

\bibitem{BKN} D.B. Boe, J. Kujawa, D.K. Nakano, Cohomology
and support varieties for Lie superalgebras, Trans. Amer. Math. Soc.
362 (2010), no. 12, 6551-6590.

\bibitem{FP1} E.M. Friedlander, B.J. Parshall, On the cohomology of
algebraic and related finite groups, Invent. Math. 74 (1983), no. 1,
85-117.

\bibitem{FP2} E.M. Friedlander, B.J. Parshall, Cohomology of
infinitesimal and discrete groups, Math. Ann. 273 (1986), no. 3,
353-374.

\bibitem{F} D.B. Fuks, Cohomology of infinite-dimensional Lie
algebras, Translated from the Russian by A. B. Sosinskii.
Contemporary Soviet Mathematics. Consultants Bureau, New York, 1986.

\bibitem{Hoch} G. Hochschild, Cohomology of restricted Lie algebras, Amer. J. Math.
76, (1954). 555-580.

\bibitem{HP} J.-S, Huang, P. Pandzic, Dirac cohomology for Lie
superalgebras. Transform. Groups 10 (2005), no. 2, 201-209.

\bibitem{Kac1} V. Kac, Lie superalgebras, Adv. in Math. 26(1977),
8-96.

\bibitem{KK} S-J. Kang, J-H. Kwon, Graded Lie superalgebras,
supertrace formula, and orbit Lie superalgebras, Proc. London Math.
Soc. (3) 81 (2000), no. 3, 675-724.

\bibitem{SF} H. Strade, R. Farnsteiner, Modular Lie algebras and
their representations, Marcel Dekker Monographs, Vol. 116(New York,
1988).

\bibitem{SZ} Y.C. Su, R.B. Zhang, Cohomology of Lie superalgebras
$\mathfrak{sl}_{m| n}$ and $\mathfrak{osp}_{2\vert 2n}$, Proc. Lond.
Math. Soc. (3) 94 (2007), no. 1, 91-136.

\bibitem{WZ} W. Wang, L. Zhao, Representations of Lie superalgebras in prime characteristic. I,
  Proc. Lond. Math. Soc. (3) 99(2009), no. 1, 145-167.


\end{thebibliography}
\end{document}